\documentclass[12pt,leqno]{article}

\def\diam{\mathop{\rm diam}}
\def\dist{\mathop{\rm dist}}
\def\osc{\mathop{\rm osc}}
\def\tr{\mathop{\rm tr}}

\begin{document}

\title{Some remarks about metric spaces}

\author{Stephen William Semmes	\\
	Rice University		\\
	Houston, Texas}

\date{}

\maketitle

\tableofcontents

\section{Basic notions and examples}
\label{basic notions and examples}
\setcounter{equation}{0}

	Of course various kinds of metric spaces arise in various
contexts and are viewed in various ways.  In this brief survey we hope
to give some modest indications of this.  In particular, we shall try
to describe some basic examples which can be of interest.

	For the record, by a \emph{metric space} we mean a nonempty
set $M$ together with a distance function $d(x, y)$, which is a
real-valued function on $M \times M$ such that $d(x, y) \ge 0$
for all $x, y \in M$, $d(x, y) = 0$ if and only if $x = y$,
$d(x, y) = d(y, x)$ for all $x, y \in M$, and
\begin{equation}
	d(x, z) \le d(x, y) + d(y, z)
\end{equation}
for all $x, y, z \in M$.  This last property is called the
\emph{triangle inequality}, and sometimes it is convenient to
allow the weaker version
\begin{equation}
	d(x, z) \le C \, (d(x, y) + d(y, z))
\end{equation}
for a nonnegative real number $C$ and all $x, y, z \in M$, which
which case $(M, d(x, y))$ is called a quasi-metric space.  Another
variant is that we may wish to allow $d(x, y) = 0$ to hold sometimes
without having $x = y$, in which case we have a semi-metric space,
or a semi-quasi-metric space, as appropriate.

	A sequence of points $\{x_j\}_{j=1}^\infty$ in a metric
space $M$ with metric $d(x, y)$ is said to converge to a point $x$
in $M$ if for every $\epsilon > 0$ there is a positive integer $L$
such that
\begin{equation}
	d(x_j, x) < \epsilon \quad\hbox{for all } j \ge L,
\end{equation}
in which case we write
\begin{equation}
	\lim_{j \to \infty} x_j = x.
\end{equation}
A sequence $\{x_j\}_{j=1}^\infty$ of points in $M$ is said to be a
\emph{Cauchy sequence} if for every $\epsilon > 0$ there is a positive
integer $L$ such that
\begin{equation}
	d(x_j, x_k) < \epsilon \quad\hbox{for all } j, k \ge L.
\end{equation}
It is easy to see that every convergent sequence is a Cauchy sequence,
and conversely a metric space in which every Cauchy sequence converges
is said to be \emph{complete}.

	A very basic example of a metric space is the real line ${\bf
R}$ with its standard metric.  Recall that if $x$ is a real number,
then the \emph{absolute value} of $x$ is denoted $|x|$ and defined to
be equal to $x$ when $x \ge 0$ and to $-x$ when $x < 0$.  One can
check that
\begin{equation}
	|x + y| \le |x| + |y|
\end{equation}
and
\begin{equation}
	|x \, y| = |x| \, |y|
\end{equation}
when $x$, $y$ are real numbers, and that the standard distance
function $|x - y|$ on ${\bf R}$ is indeed a metric.

	Let $n$ be a positive integer, and let ${\bf R}^n$ denote the
real vector space of $n$-tuples of real numbers.  Thus elements $x$ of
${\bf R}^n$ are of the form $(x_1, \ldots, x_n)$, where the $n$
coordinates $x_j$, $1 \le j \le n$, are real numbers.  If $x = (x_1,
\ldots, x_n)$, $y = (y_1, \ldots, y_n)$ are two elements of ${\bf
R}^n$ and $r$ is a real number, then the sum $x + y$ and scalar
product $r \, x$ are defined coordinatewise in the usual manner, by
\begin{equation}
	x + y = (x_1 + y_1, \ldots, x_n + y_n)
\end{equation}
and
\begin{equation}
	r \, x = (r \, x_1, \ldots, r \, x_n).
\end{equation}

	If $x$ is an element of ${\bf R}^n$, then the standard
Euclidean norm of $x$ is denoted $|x|$ and defined by
\begin{equation}
	|x| = \biggl(\sum_{j=1}^n x_j^2 \biggr)^{1/2}.
\end{equation}
One can show that
\begin{equation}
	|x + y| \le |x| + |y|
\end{equation}
holds for all $x, y \in {\bf R}^n$, and we shall come back to this
in a moment, and clearly we also have that
\begin{equation}
	|r \, x| = |r| \, |x|
\end{equation}
for all $r \in {\bf R}$ and $x \in {\bf R}^n$, which is to say that
the norm of a scalar product of a real number and an element of ${\bf
R}^n$ is equal to the product of the absolute value of the real number
and the norm of the element of ${\bf R}^n$.  Using these properties,
one can check that the standard Euclidean distance $|x - y|$ on
${\bf R}^n$ is indeed a metric.

	More generally, a \emph{norm} on ${\bf R}^n$ is a real-valued
function $N(x)$ such that $N(x) \ge 0$ for all $x \in {\bf R}^n$,
$N(x) = 0$ if and only if $x = 0$, 
\begin{equation}
	N(r \, x) = |r| \, N(x)
\end{equation}
for all $r \in {\bf R}$ and $x \in {\bf R}^n$, and
\begin{equation}
	N(x + y) \le N(x) + N(y)
\end{equation}
for all $x, y \in {\bf R}^n$.  If $N(x)$ is a norm on ${\bf R}^n$,
then
\begin{equation}
	d(x, y) = N(x - y)
\end{equation}
defines a metric on ${\bf R}^n$.  As for metrics, one can weaken the
triangle inequality or relax the condition that $N(x) = 0$ implies
$x = 0$ to get quasi-norms, semi-norms, and semi-quasi-norms.

	Recall that a subset $E$ of ${\bf R}^n$ is said to be
\emph{convex} if
\begin{equation}
	t \, x + (1 - t) \, y \in E
\end{equation}
whenever $x$, $y$ are elements of $E$ and $t$ is a real number such
that $0 < t < 1$.  A real-valued function $f(x)$ on ${\bf R}^n$ is
said to be convex if and only if
\begin{equation}
	f(t \, x + (1 - t) \, y) \le t \, f(x) + (1 - t) \, f(y)
\end{equation}
for all $x, y \in {\bf R}^n$ and $t \in {\bf R}^n$ with $0 < t < 1$.
If $N(x)$ is a real-valued function on ${\bf R}^n$ which is assumed to
satisfy the conditions of a norm except for the triangle inequality,
then one can check that the triangle inequality, the convexity of the
closed unit ball
\begin{equation}
	\{x \in {\bf R}^n : N(x) \le 1\},
\end{equation}
and the convexity of $N(x)$ as a function on ${\bf R}^n$, are all
equivalent.

	For example, if $p$ is a real number such that $1 \le p <
\infty$, then define $|x|_p$ for $x \in {\bf R}^n$ by
\begin{equation}
	|x|_p = \biggl(\sum_{j=1}^n |x|^p \biggr)^{1/p},
\end{equation}
which is the same as the standard norm $|x|$ when $p = 2$.  For
$p = \infty$ let us set
\begin{equation}
	|x|_\infty = \max \{|x_j| : 1 \le j \le n\}.
\end{equation}
One can check that these define norms on ${\bf R}^n$, using the
convexity of the function $|r|^p$ on ${\bf R}$ when $1 < p < \infty$
to check that the closed unit ball of $|x|_p$ is convex and hence that
the triangle inequality holds when $1 < p < \infty$.

	Let us now consider a class of metric spaces along the lines
of Cantor sets.  For this we assume that we are given a sequence
$\{F_j\}_{j=1}^\infty$ of nonempty finite sets.  We also assume that
$\{\rho_j\}_{j=1}^\infty$ is a monotone decreasing sequence of
positive real numbers which converges to $0$.

	For our space $M$ we take the Cartesian product of the
$F_j$'s, so that en element $x$ of $M$ is a sequence
$\{x_j\}_{j=1}^\infty$ such that $x_j \in F_j$ for all $j$.
We define a distance function $d(x, y)$ on $M$ by setting
$d(x, y) = 0$ when $x = y$, and
\begin{equation}
	d(x, y) = \rho_j
\end{equation}
when $x_j \ne y_j$ and $x_i = y_i$ for all $i < j$.  One can check
that this does indeed define a metric space, and in fact $d(x,y)$
is an \emph{ultrametric}, which is to say that
\begin{equation}
	d(x, z) \le \max(d(x, y), d(y, z))
\end{equation}
for all $x, y, z \in M$.

	The classical Cantor set is the subset of the unit interval
$[0, 1]$ in the real line obtained by first removing the open
subinterval $(1/3, 2/3)$, then removing the the open middle thirds
of the two closed intervals which remain, and so on.  Alternatively,
the classical Cantor set can be described as the set of real numbers
$t$ such that $0 \le t \le 1$ and $t$ has ab expansion base $3$
whose coefficients are are either $0$ or $2$.  This set, equipped
with the standard Euclidean metric, is very similar to the general
situation just described with each $F_j$ having two elements and
with $\rho_j = 2^{-j}$ for all $j$, although the metrics are not
quite the same.

	In general, if $(M, d(x, y))$ is a metric space and $E$ is a
nonempty subset of $M$, then $E$ can be considered as a metric space
in its own right, using the restriction of the metric $d(x, y)$ from
$M$ to $E$.  Sometimes there may be another metric on $E$ which is
similar to the one inherited from the larger space $M$, and which has
other nice properties, as in the case of Cantor sets just described.
Another basic instance of this occurs with arcs in Euclidean spaces
which are ``snowflakes'', and which are similar to taking the unit
interval $[0, 1]$ in the real line with the metric $|x - y|^a$
for some real number $a$, $0 < a < 1$, or other functions of the
standard distance on $[0, 1]$.

	A nonempty subset $E$ of a metric space $(M, d(x,y))$ is said
to be \emph{bounded} if the real numbers $d(x, y)$, $x, y \in E$, are
bounded, in which case the \emph{diameter} of $E$ is denoted $\diam E$
and defined by
\begin{equation}
	\diam E = \sup \{d(x, y) : x, y \in E \}.
\end{equation}
A stronger condition is that $E$ be \emph{totally bounded}, which
means that for each $\epsilon > 0$ there is a finite collection $A_1,
\ldots, A_k$ of subsets of $E$ such that
\begin{equation}
	E \subseteq \bigcup_{j=1}^k A_j
\end{equation}
and
\begin{equation}
	d(x, y) < \epsilon \quad x, y \in A_j,
\end{equation}
$j = 1, \ldots, k$.  A basic feature of Euclidean spaces is that
bounded subsets are totally bounded, and the generalized Cantor
sets described before are totally bounded.

	A metric space $(M, d(x, y))$ is \emph{compact} if it is
complete, so that every Cauchy sequence converges, and totally
bounded.  This is equivalent to the standard definitions in terms of
open coverings or the existence of limit points.  A closed and bounded
subset of ${\bf R}^n$ is compact, and the generalized Cantor sets
described earlier are compact.

	Another way that metric spaces arise is to start with a
connected smooth $n$-dimensional manifold $M$, which is basically a
space which looks locally like $n$-dimensional Euclidean space.  At
each point $p$ in $M$ one has an $n$-dimensional tangent space
$T_p(M)$, which looks like ${\bf R}^n$ as a vector space, and on which
one can put a norm.  If at each point $p$ in $M$ one can identify
$T_p(M)$ with ${\bf R}^n$ with its standard norm, then the space is
Riemannian, and with general norms the space is of Riemann--Finsler
type.

	In this type of situation, the length of a nice path in $M$
can be defined by integrating the infinitesimal lengths determined by
the norms on the tangent spaces.  The distance between two points is
defined to be the infimum of the lengths of the paths connecting the
two points.  It is easy to see that this does indeed define a metric,
with the triangle inequality being a consequence of the way that the
distance is defined.

	A basic example of this is the $n$-dimensional sphere
${\bf S}^n$ in ${\bf R}^{n+1}$, defined by
\begin{equation}
	{\bf S}^n = \{x \in {\bf R}^{n+1} : |x| = 1\}.
\end{equation}
The tangent space of ${\bf S}^n$ can be identified with the
$n$-dimensional linear subspace of ${\bf R}^{n+1}$ of vectors $x$
which are orthogonal to $p$ as a vector itself.  This leads to a
Euclidean norm on the tangent spaces, inherited from the one on ${\bf
R}^{n+1}$.

	One can also define distances through paths in other
situations.  In a nonempty connected graph, every pair of points is
connected by a path, the length of a path can be defined as the number
of edges traversed, and the distance between two points can be defined
as the length of the shortest path between the two points.  In many
standard fractals, like the Sierpinski gasket and carpet, there are
a lot of paths of finite length between arbitrary elements of the
fractal, and the infimum of the lengths of these paths defines a
metric on the fractal.

	Let $(M, d(x, y))$ be a metric space.  If $A$, $B$ are
nonempty subsets of $M$ and $t$ is a positive real number, then let us
say that $A$, $B$ are ``$t$-close'' if for every $x \in A$ there is a
$y \in B$ such that $d(x, y) < t$, and if for every $y \in B$ there is
an $x \in A$ such that $d(x, y) < t$.  By definition, this relation is
symmetric in $A$ and $B$.

	A subset $E$ of $M$ is said to be \emph{closed} if for
every sequence $\{z_j\}_{j=1}^\infty$ of points in $E$ which
converges to a point $z$ in $M$, we have that $z \in E$.
Let us write $\mathcal{S}(M)$ for the set of nonempty closed
and bounded subsets of $M$.  If $A$, $B$ are two elements of
$\mathcal{S}(M)$, then the \emph{Hausdorff distance} between
$A$ and $B$ is defined to be the infimum of the set of positive
real numbers $t$ such that $A$, $B$ are $t$-close.

	Because $A$ and $B$ are bounded subsets of $M$, there are
positive real numbers $t$ such that $A$, $B$ are $t$-close.  The
restriction to closed sets ensures that if $A$, $B$ are $t$-close for
all $t > 0$, then $A = B$.  If $A_1$, $A_2$, $A_3$ are nonempty
subsets of $M$ and $t_1$, $t_2$ are positive real numbers such that
$A_1$, $A_2$ are $t_1$-close and $A_2$, $A_3$ are $t_2$-close, then
one can check that $A_1$, $A_3$ are $(t_1 + t_2)$-close, and this
implies the triangle inequality for the Hausdorff distance.

	From this it follows easily that $\mathcal{S}(M)$, equipped
with the Hausdorff distance $D(A, B)$, is indeed a metric space.  A
basic result states that if $M$ is compact, then $\mathcal{S}(M)$ is
compact too.  This is not too difficult to show.

	Now let us consider a situation with a lot of deep
mathematical structure which has been much-studied, involving
interplay between algebra, analysis, and geometry.

	Fix a positive integer $n$, which we may as well take to be at
least $2$.  Let $\mathcal{M}_n^+$ denote the set of $n \times n$ real
symmetric matrices which are positive-definite and have determinant
equal to $1$.  One can think of this as a smooth hypersurface in the
vector space of $n \times n$ real symmetric matrices, and it is a
smooth manifold in particular.

	We can start with a Riemannian view of this space.  For each
$H \in \mathcal{M}_n^+$, we can identify the tangent space of
$\mathcal{M}_n^+$ at $H$ with the vector space of $n \times n$ real
symmetric matrices $A$ such that the trace of $H^{-1} \, A$ is equal
to $0$.  Thus for each such $A$, we get a one-parameter family of
perturbations of $H$ by taking $H + t \, A$, where $t$ is a real
number with small absolute value, so that $H + t \, A$ is still
positive definite, and to first order in $t$ these perturbations also
have determinant equal to $1$.

	Conversely, to first order in $t$, each smooth perturbation of
$H$ in $\mathcal{M}_n^+$ is of this form.  Now, for each $A$ of this
type, we define its norm as an element of the tangent space to
$\mathcal{M}_n^+$ to be
\begin{equation}
	\biggl(\tr H^{-1} \, A \, H^{-1} \, A\biggl)^{1/2}.
\end{equation}
Here $\tr B$ denotes the trace of a square matrix $A$, and we are
using ordinary matrix multiplication in this expression.

	When $H$ is the identity matrix, this reduces to the square
root of $\tr A^2$, which is a kind of Euclidean norm of $A$.  For
general $H$'s, we adapt the norm to $H$.  It is still basically
a Euclidean norm, so that we are in the Riemannian case.

	Let us consider the transformation on $\mathcal{M}_n^+$
defined by $H \mapsto H^{-1}$.  If $H + t \, A$ is a basic
first-order deformation of $H$, as before, then that is transformed
to $(H + t \, A)^{-1}$, which is the same as
\begin{equation}
	H^{-1} - t \, H^{-1} \, A \, H^{-1}
\end{equation}
to first order in $t$.  Thus $A$ as a tangent vector at $H$
corresponds to $- H^{-1} \, A \, H^{-1}$ as a tangent vector at
$H^{-1}$ under the mapping $H \mapsto H^{-1}$, and it is easy to see
that the norm of $A$ as a tangent vector at $H$ is equal to the norm
of $-H^{-1} \, A \, H^{-1}$ as a tangent vector at $H^{-1}$.

	Now let $T$ be an $n \times n$ matrix with determinant equal
to $1$, and let $T^*$ denote its transpose, which is also an
invertible $n \times n$ matrix.  Associated to $T$ is the mapping
\begin{equation}
	H \mapsto T \, H \, T^*,
\end{equation}
and if $A$ corresponds to a tangent vector at $H$ as before, then
$T \, A \, T^*$ is the tangent vector at $T \, H \, T^*$ induced
by our mapping.  Again one can check that the norm of $A$ as a
tangent vector at $H$ is the same as the norm of $T \, A \, T^*$
as a tangent vector at $T \, H \, T^*$.

	Assume further that the entries of $T$ are integers.  This
implies that the inverse of $T$ also has integer entries, by Cramer's
rule.  The product of two such matrices has the same property, and
indeed this defines a nice discrete group of matrices.

	This discrete group acts on $\mathcal{M}_n^+$, and we can pass
to the corresponding quotient space.  That is, we now identify two
elements $H_1$, $H_2$ in $\mathcal{M}_n^+$ if there is an integer
matrix $T$ as above so that
\begin{equation}
	H_2 = T \, H_1 \, T^*.
\end{equation}
This relation between $H_1, H_2 \in \mathcal{M}_n^+$ is indeed an
equivalence relation, so that we get a nice quotient.

	Because of the discreteness of the group of integer matrices
with determinant $1$, the quotient space is still a smooth manifold,
since it looks locally like $\mathcal{M}_n^+$.  The norm on the
tangent spaces still makes sense as well, because the transformations
defining the equivalence relation preserves these norms, as we have
seen.  Thus this quotient of $\mathcal{M}_n^+$ is still a nice smooth
connected Riemannian manifold.

	In some situations like this the quotient space turns out to
be compact.  In this case the quotient space is not compact, but it
does have finite volume.  More precisely, the notion of volume at the
level of the tangent spaces is determined by the norm, and preserved
in the present circumstances by the transformations used in the
equivalence relation, and the volume of the quotient can be obtained
by integrating the infinitesimal volumes.

	The group of $n \times n$ matrices with integer entries and
determinant $1$ is a very interesting special case of discrete groups
more generally.  Suppose now that $\Gamma$ is a group and that $A$
is a finite symmetric subset of $\Gamma$ which generates $\Gamma$,
so that $\alpha^{-1} \in A$ when $\alpha \in A$ and every element
of $\Gamma$ can be expressed as a finite product of elements of $A$,
with the identity element automatically corresponding to an empty
product.  This leads to the associated \emph{Cayley graph}, in which
two elements $\gamma_1$, $\gamma_2$ of $\Gamma$ are considered to be
adjacent if $\gamma_2$ can be expressed as $\gamma_1 \, \alpha$
for some $\alpha \in A$, and to a distance function on $\Gamma$
which is invariant under left-multiplication in the group.

\section{Aspects of analysis}
\label{aspects of analysis}
\setcounter{equation}{0}

	Let $(M, d(x, y))$ be a metric space, let $f(x)$ be a real-valued
function on $M$, and let $C$ be a nonnegative real number.  We say that
$f$ is \emph{$C$-Lipschitz} if
\begin{equation}
	|f(x) - f(y)| \le C \, d(x, y)
\end{equation}
for all $x, y \in M$.  This is equivalent to saying that
\begin{equation}
	f(x) \le f(y) + C \, d(x, y)
\end{equation}
for all $x, y \in M$.  Notice that a function is $0$-Lipschitz if and
only if it is constant.

	For instance, if $p$ is a point in $M$, then
$f_p(x) = d(x, p)$ is $1$-Lipschitz, because
\begin{equation}
	d(x, p) \le d(x, y) + d(y, p)
\end{equation}
by the triangle inequality.  More generally, if $A$ is a nonempty
subset of $M$, and if $x$ is a point in $M$, then the distance from
$x$ to $A$ is denoted $\dist(x, A)$ and defined by
\begin{equation}
	\dist(x, A) = \inf \{d(x, y) : y \in A\},
\end{equation}
and one can check that $\dist(x, A)$ is a $1$-Lipschitz function on
$M$.  If $f_1$, $f_2$ are two real-valued functions on $M$ which are
$C_1$, $C_2$-Lipschitz, respectively, and if $\alpha_1$, $\alpha_2$
are real numbers, then $\max(f_1, f_2)$, $\min(f_1, f_2)$ are
$C$-Lipschitz with $C = \max (C_1, C_2)$, and $\alpha_1 \, f_1 +
\alpha_2 \, f_2$ is $C$-Lipschitz with $C = |\alpha_1| \, C_1 +
|\alpha_2| \, C_2$.

	Now let $C$ be a nonnegative real number and let $s$ be a
positive real number.  A real-valued function $f(x)$ on $M$ is said to
be \emph{$C$-Lipschitz of order $s$} if
\begin{equation}
	|f(x) - f(y)| \le C \, d(x, y)^s
\end{equation}
for all $x, y \in M$, which is again equivalent to
\begin{equation}
	f(x) \le f(y) + C \, d(x, y)^s
\end{equation}
for all $x, y \in M$.  As before, $f(x)$ is $0$-Lipschitz of order $s$
if and only if $f(x)$ is constant on $M$.

	When $0 < s < 1$, one can check that $d(x, y)^s$ is also a
metric on $M$, which defines the same topology on $M$ in fact.  The
main point in this regard is that the triangle inequality continues to
hold, which follows from the observation that
\begin{equation}
	(\alpha + \beta)^s \le \alpha^s + \beta^s
\end{equation}
for all nonnegative real numbers $\alpha$, $\beta$.  A real-valued
function $f(x)$ on $M$ is $C$-Lipschitz of order $s$ with respect to
the metric $d(x, y)$ if and only if $f(x)$ is $C$-Lipschitz of order
$1$ with respect to $d(x, y)^s$, and as a result when $0 < s < 1$ one
has the same statements for Lipschitz functions of order $s$ as for
ordinary Lipschitz functions.

	When $s > 1$ the triangle inequality for $d(x, y)^s$ does not
work in general, but we do have that
\begin{equation}
	d(x, z)^s \le 2^{s-1} \, (d(x, y)^s + d(y, z)^s)
\end{equation}
for all $x, y, z \in M$, because
\begin{equation}
	(\alpha + \beta)^s \le 2^{s-1} \, (\alpha^s + \beta^s)
\end{equation}
for all nonnegative real numbers $\alpha$, $\beta$.  Some of the usual
properties of Lipschitz functions carry over to Lipschitz functions of
order $s$, perhaps with appropriate modification, but it may be that
the only Lipschitz functions of order $s$ when $s > 1$ are constant.
This is the case on Euclidean spaces with their standard metrics,
and indeed a Lipschitz function of order $s > 1$ has first derivatives
equal to $0$ everywhere.

	In harmonic analysis one considers a variety of classes of
functions with different kinds of restrictions on size, oscillations,
regularity, and so on, and these Lipschitz classes are fundamental
examples.  In particular, it can be quite useful to have the parameter
$s$ available to adjust to the given circumstances.  There are also
other ways of introducing parameters to get interesting classes of
functions and measurements of their behavior.

	If $M$ is the usual $n$-dimensional Euclidean space ${\bf
R}^n$, with its standard metric, then one has the extra structure of
translations, rotations, and dilations.  If $f(x)$ is a real-valued
function on ${\bf R}^n$ which is $C$-Lipschitz of order $s$, $f(x -
u)$ is also $C$-Lipschitz of order $s$ for each $u \in {\bf R}^n$,
$f(\Theta(x))$ is $C$-Lipschitz of order $s$ for each rotation
$\Theta$ on ${\bf R}^n$, and $f(r^{-1} x)$ is $(C \, r^s)$-Lipschitz
of order $s$ for each $r > 0$.  In effect, on general metric spaces we
can consider classes of functions and measurements of their behavior
which have analogous features, even if there are not exactly
translations, rotations, and dilations.

	A basic notion is to consider various scales and locations
somewhat independently.  In this regard, if $f(x)$ is a real-valued
function on $M$, $x$ is an element of $M$, and $t$ is a positive real
number, put
\begin{equation}
	\osc(x, t) = \sup \{|f(y) - f(x)| : y \in M, d(y, x) \le t\}.
\end{equation}
We implicitly assume here that $f(y)$ remains bounded on bounded
subsets of $M$, so that this quantity is finite.  For instance, $f$ is
$C$-Lipschitz of order $s$ if and only if
\begin{equation}
	t^{-s} \osc(x, t) \le C
\end{equation}
for all $x \in M$ and $t > 0$.  

	Let us pause a moment and notice that
\begin{equation}
	\osc(w, r) \le \osc(x, t)
\end{equation}
when $d(w, x) + r \le t$.  Thus, 
\begin{equation}
	r^{-s} \, \osc(w, r) \le 2^s \, t^{-s} \, \osc(x, t)
\end{equation}
when $d(x, w) + r \le t$ and $r \ge t/2$.  This is a kind of
``robustness'' property of these measurements of local oscillation of
a function $f$ on $M$.  In particular, to sample the behavior of $f$
at essentially all locations and scales, it is practically enough to
look at a reasonably-nice and discrete family of locations and scales.
For example, one might restrict one's attention to radii $t$ which are
integer powers of $2$, and for a specific choice of $t$ use a
collection of points in $M$ which cover suitably the various locations
at that scale.

	Instead of simply taking a supremum of some measurements of
local oscillation like this, one can consider various sums of discrete
samples of this sort.  This leads to a number of classes of functions
and measurements of their behavior.  One can adjust this further
by taking into account the relation of some location and scale
to some kind of boundaries, or singularities, or concentrations,
and so on.

	Of course one might also use some kind of measurement of sizes
of subsets of $M$.  This could entail diameters, volumes, or
measurements of capacity.  There are also many kinds of local
measurements of oscillation or size that one can consider.  As an
extension of just taking suprema, one can take various local averages,
or averages of powers of other quantities.  Of course one can still
bring in powers of the radius as before.

	Even if one starts with measurements of localized behavior
which are not so robust in the manner described before, one can
transform them into more robust versions by taking localized suprema
or averages or whatever afterwards.  Frequently the kind of overall
aggregations employed have this kind of robustness included in effect,
and one can make some sort of rearrangement to put this in starker
relief.  Let us also note that one often has local measurements which
can be quite different on their own, but in some overall aggregation
lead to equivalent classes of functions and similar measurements of
their behavior.

	There are various moments, differences, and higher-order
oscillations that can be interesting.  As a basic version of this, one
can consider oscillations of $f(x)$ in terms of deviations from
something like a polynomial of fixed positive degree, rather than
simply oscillations from being constant, as with $\osc(x, t)$.  This
can be measured in a number of ways.

	However, for these kinds of higher-order oscillations,
additional structure of the metric space is relevant.  On Euclidean
spaces, or subsets of Euclidean spaces, one can use ordinary
polynomials, for instance.  This carries over to the much-studied
setting of nilpotent Lie groups equipped with a family of dilations,
where one has polynomials as in the Euclidean case, with the
degrees of the polynomials defined in a different way using the
dilations.

	These themes are closely related to having some kind of
derivatives around.  Just as there are various ways to measure the
size of a function, one can get various measurements of oscillations
looking at measurements of sizes of derivatives.  It can also be
interesting to have scales involved in a more active manner, and in
any case there are numerous versions of ideas along these lines
that one can consider.

\section{Sub-Riemannian geometry}
\label{sub-riemannian geometry}
\setcounter{equation}{0}

	Let $n$ be a positive integer, and consider $(2n+1)$-dimensional
Euclidean space, which we shall think of as
\begin{equation}
\label{2n+1 dimensional space}
	{\bf R}^n \times {\bf R}^n \times {\bf R}.
\end{equation}
The case where $n = 1$, so that this is basically just ${\bf R}^3$,
is already quite interesting for us.  We shall also be interested in
\begin{equation}
\label{2n dimensional space}
	{\bf R}^n \times {\bf R}^n
\end{equation}
and the obvious coordinate projection from the former to the latter
given by
\begin{equation}
	(x, y, s) \mapsto (x, y).
\end{equation}

	Define a smooth $1$-form $\alpha$ on (\ref{2n+1 dimensional
space}) by
\begin{equation}
	\alpha = ds - \sum_{j = 1}^n y_j \, dx_j,
\end{equation}
using the usual coordinates $(x, y, s)$ on (\ref{2n+1 dimensional
space}).  Thus if $p = (x, y, s)$ is a point in (\ref{2n+1 dimensional
space}), then $\alpha$ at $p$, which we denote $\alpha_p$, is a linear
functional on the tangent space to (\ref{2n+1 dimensional space}) at
$p$.  Of course we can identify the tangent space to (\ref{2n+1
dimensional space}) at $p$ with ${\bf R}^n \times {\bf R}^n \times
{\bf R}$ itself in the usual way, and if $V = (v, w, u)$ is a tangent
vector to (\ref{2n+1 dimensional space}) at $p$ represented using this
identification, then
\begin{equation}
	\alpha_p(V) = u - \sum_{j=1}^n y_j \, v_j.
\end{equation}

	We can take the derivative $d \alpha$ of $\alpha$ in the sense
of exterior differential calculus to get a $2$-form on (\ref{2n+1
dimensional space}).  Namely,
\begin{equation}
	d \alpha = \sum_{j=1}^n dx_j \wedge dy_j.
\end{equation}
If one takes the wedge product of $\alpha$ with $n$ copies of $d \alpha$,
then one gets
\begin{equation}
	(n!) \, ds \wedge dx_1 \wedge dy_1 \cdots dx_n \wedge dy_n,
\end{equation}
which is a nonzero $(2n+1)$-form on (\ref{2n+1 dimensional space})
that is $n$ factorial times the standard volume form.

	If $p = (x, y, s)$ is a point in (\ref{2n+1 dimensional space}),
then we get a special linear subspace $H_p$ of the tangent space to
(\ref{2n+1 dimensional space}) at $p$ which is the kernel of $\alpha_p$.
In other words, $H_p$ consists of the tangent vectors $V = (v, w, u)$
such that
\begin{equation}
	u - \sum_{j=1}^n y_j \, v_j = 0.
\end{equation}
Thus $H_p$ has dimension $2n$ at every point $p$.

	If $f$ is a smooth real-valued function on (\ref{2n+1
dimensional space}) which is never equal to $0$, then $f \, \alpha$
is also a $1$-form on (\ref{2n+1 dimensional space}) which vanishes
exactly at the tangent vectors in $H_p$ at a point $p$.  By the usual
rules of exterior differential calculus,
\begin{equation}
	d(f \, \alpha) = df \wedge \alpha + f \, d \alpha.
\end{equation}
One can check that the wedge product of $f \alpha$ with $n$ copies
of $d(f \alpha)$ is the same as $f^{n+1}$ times the wedge product of
$\alpha$ with $n$ copies of $d \alpha$.

	A nice feature of $d \alpha$ is that it can be viewed as the
pull-back of a $2$-form from (\ref{2n dimensional space}).  Basically
this simply means that $d \alpha$ does not contain any $ds$'s and the
coefficients depend only on $x$, $y$.  However, $\alpha$ contains $ds$
as an important term, and $\alpha$ is not the pull-back of a form
from (\ref{2n dimensional space}).

	For each point $p = (x, y, s)$ in (\ref{2n+1 dimensional
space}), let $F_p$ denote the $1$-dimensional linear subspace of the
tangent space to (\ref{2n+1 dimensional space}) at $p$ consisting of
vectors of the form $(0, 0, t)$, $t \in {\bf R}$.  This is the same as
the subspace of tangent vectors to the fiber of the mapping from
(\ref{2n+1 dimensional space}) to (\ref{2n dimensional space}) through
$p$, which is also the same as the set of tangent vectors in the
kernel of the aforementioned mapping.  Notice that $F_p$ and $H_p$ are
complementary subspaces of the tangent space to (\ref{2n+1 dimensional
space}), which is to say that every vector $V = (v, w, u)$ in the
tangent space at $p$ can be written in a unique way as a sum of
vectors in $F_p$ and $H_p$, namely
\begin{equation}
	V = \biggl(0, 0, u - \sum_{j=1}^n y_j \, v_j \biggr)
		+ \biggl(v, w, \sum_{j=1}^n y_j \, v_j \biggr).
\end{equation}

	Let $I$ be an interval in the real line with positive length,
and which may be unbounded.  Suppose that $\gamma(t)$ is a continuous
function defined for $t$ in $I$ and with values in (\ref{2n+1
dimensional space}), so that $\gamma(t)$ defines a continuous path in
(\ref{2n+1 dimensional space}).  Let us assume that $\gamma(t)$ is
continuously-differentiable on $I$, so that the derivative
\begin{equation}
	\dot{\gamma}(t) = \frac{d}{dt} \gamma(t)
\end{equation}
exists for each $t \in I$ and is continuous on $I$.  More precisely,
if $t$ is in the interior of $I$, then the derivative is taken in the
usual sense, while if $t$ is an element of $I$ which is also an
endpoint of $I$, then one uses a one-sided derivative.

	We say that this path $\gamma(t)$ is \emph{horizontal} if
\begin{equation}
	\dot{\gamma}(t) \in H_{\gamma(t)}
\end{equation}
for each $t \in I$.  This is equivalent to asking that
\begin{equation}
	\alpha(\dot{\gamma}) = 0
\end{equation}
along the curve, which is to say that
\begin{equation}
	\alpha_{\gamma(t)}(\dot{\gamma}(t)) = 0
\end{equation}
for all $t \in I$.  If we write $\gamma(t)$ more explicitly as
$(x(t), y(t), s(t))$, then this becomes
\begin{equation}
	\dot{s}(t) = \sum_{j=1}^n y_j(t) \, \dot{x}_j(t)
\end{equation}
for all $t \in I$.

	Let $\gamma_0(t)$ be the projection of $\gamma(t)$ into
(\ref{2n dimensional space}), which means that $\gamma_0(t) = (x(t),
y(t))$.  If $\gamma(t)$ is a horizontal curve in (\ref{2n+1
dimensional space}), then $\gamma(t)$ is uniquely determined by its
projection $\gamma_0(t)$ and the value of $\gamma(t)$ at a single
point $t = t_0 \in I$.  Indeed, if $\gamma(t) = (x(t), y(t), s(t))$ is
horizontal, then $\dot{s}(t)$ can be expressed in terms of
$\dot{x}(t)$ and $y(t)$ as before, and $s(t)$ is determined by this
and the value of $s(t)$ at one point $t_0$.

	Conversely, suppose that we are given a
continuously-differentiable curve $\gamma_0(t)$ in (\ref{2n
dimensional space}) defined for $t \in I$.  Also let $t_0$ be an
element of $I$, and suppose that $s_0$ is some real number.  Then
there is a continuously-differentiable curve $\gamma(t) = (x(t), y(t),
s(t))$, $t \in I$, in (\ref{2n+1 dimensional space}) whose projection
into (\ref{2n dimensional space}) is equal to $\gamma_0(t)$ and such
that $s(t_0) = s_0$.  Namely, one can compute what $\dot{s}(t)$ should
be in terms of $\dot{x}(t)$ and $y(t)$, and integrate that using also
$s(t_0) = s_0$ to get $s(t)$ for all $t \in I$.

	It is not too difficult to show that each pair of points $p$,
$q$ in (\ref{2n+1 dimensional space}) can be connected by a
continuously-differentiable curve which is horizontal.  One can start
by taking the prjections of these two points to get two points in
(\ref{2n dimensional space}) which can be connected by all sorts of
curves.  These curves can be lifted to horizontal curves in (\ref{2n+1
dimensional space}) which begin at $p$, as in the preceding paragraph.
In general a lifted path like this may not go to $q$, and one can
check that there are plenty of choices of paths in (\ref{2n
dimensional space}) for which the lifting will have this property.

	This leads to a very interesting kind of geometry in
(\ref{2n+1 dimensional space}).  Namely, one defines the distance
between two points $p$, $q$ to be the infimum of the lengths of the
horizontal paths joining $p$ to $q$.  Of course the standard Euclidean
metric on (\ref{2n+1 dimensional space}) can be defined by minimizing
the lengths of all paths from $p$ to $q$.  The restriction to horizontal
paths makes the distance increase, although one can check that the
resulting metric is still compatible with the usual topology on
(\ref{2n+1 dimensional space}).

\section{Hyperbolic groups}
\label{hyperbolic groups}
\setcounter{equation}{0}

	Let $\Gamma$ be a group, and let $F$ be a finite set of
elements of $\Gamma$.  By a \emph{word} over $F$ we mean a formal
product of elements of $F$ and their inverses.  Every word over $F$
determines an element of the group $\Gamma$, simply using the group
operations.  The ``empty word'' is considered a word over $F$, which
corresponds to the identity element of $\Gamma$.

	If $z$ is a word over $F$, then the length of $z$ is denoted
$L(z)$ and is the number of elements of $F$ such that they or their
inverses are used in $z$, counting multiplicities.  A word $z$ is said
to be irreducible if it does not contain an $\alpha \in F$ next to
$\alpha^{-1}$, i.e., so that all obvious cancellations have been made.
If a word $z$ over $F$ corresponds to the identity element of $\Gamma$,
then $z$ is said to be trivial.

	A finite subset $F$ of a group $\Gamma$ is a set of
\emph{generators} of $\Gamma$ if every element of $\Gamma$ corresponds
to a word over $F$.  A group is said to be \emph{finitely-generated}
if it has a finite set of generators.  Let us make the convention that
a generating set $F$ of a group $\Gamma$ should not contain the
identity element of $\Gamma$.

	Suppose that $\Gamma$ is a group and that $F$ is a finite set
of generators of $\Gamma$.  The \emph{Cayley graph} associated to
$\Gamma$ and $F$ is the graph consisting of the elements of $\Gamma$
as vertices with the provision that $\gamma_1$, $\gamma_2$ in $\Gamma$
are adjacent if $\gamma_2 = \gamma_1 \, \alpha$, where $\alpha$ is an
element of $F$ or its inverse.  Thus this relation is symmetric in
$\gamma_1$ and $\gamma_2$.

	A finite sequence $\theta_0, \theta_2, \ldots, \theta_k$ of
elements of $\Gamma$ is said to define a \emph{path} if $\theta_j$,
$\theta_{j+1}$ are adjacent in the Cayley graph for each $j$, $0 \le j
\le k-1$.  The \emph{length} of this path is defined to be $k$.  We
include the degenerate case where $k = 0$, so that a single element of
$\Gamma$ is viewed as a path of length $0$.

	If $\phi$, $\psi$ are elements of $\Gamma$, then the
\emph{distance} between $\phi$ and $\psi$ is defined to be the
shortest length of a path that connects $\phi$ to $\psi$.  In
particular, note that for any two elements $\phi$, $\psi$ in $\Gamma$
there is a path which starts at $\phi$ and ends at $\psi$.  To see
this, one can write $\psi$ as $\phi \, \beta$ for some $\beta$ in
$\Gamma$, and then use the assumption that $\Gamma$ is generated by
$F$ to obtain a path from $\phi$ to $\psi$ one step at a time.

	These definitions are invariant under \emph{left} translations
in $\Gamma$.  In other words, if $\delta$ is any fixed element of
$\Gamma$, then the tranformation $\gamma \mapsto \delta \, \gamma$ on
$\Gamma$ defines an automorphism of the Cayley graph, and it also
preserves distances between elements of $\Gamma$.  This follows
from the definitions, since the Cayley graph was defined in terms of
right-multiplication by generators and their inverses.

	A basic fact is that this definition of distance does not
depend too strongly on the choice of generating set $F$, in the sense
that if one has another finite generating set, then the two distance
functions associated to these generating sets are each bounded by a
constant multiple of the other.  This is not difficult to check, by
expressing each generator in one set as a finite word over the other
set of generators.  There are only a finite number of these expressions,
so that their maximal length is a finite number.

	Let us continue with the assumption that we have a fixed
generating set $F$ for the group $\Gamma$.  Suppose that $R$ is a
finite set of words over $F$.  We say that $R$ is a set of
\emph{relations} for $\Gamma$ if every element of $R$ is a trivial
word.  The inverses of elements of $R$ are also then trivial words, as
well as conjugates of elements of $R$.  That is, if $r$ is an element
of $R$ and $u$ is any word over $F$, then $u \, r \, u^{-1}$ is the
conjugate of $r$ by $u$, and it is a trivial word since $r$ is.
Products of conjugates of elements of $R$ and their inverses are
trivial words too, as well as words obtained from these through
cancellations, i.e., by cancelling $\alpha \, \alpha^{-1}$ and
$\alpha^{-1} \, \alpha$ whenever $\alpha$ is an element of $F$.  The
combination of $F$ and a set $R$ of relations defines a
\emph{presentation} of $\Gamma$ if every word over $F$ which
corresponds to the identity element of $\Gamma$ can be obtained in
this manner.  The empty word is viewed as being equal to the empty
product of relations, so that it is automatically included.  A group
$\Gamma$ is said to be \emph{finitely-presented} if there is a
presentation with a finite set of generators and a finite set of
relations.  For instance, if $\Gamma$ is the free group with
generators in $F$, then one can take $R$ to be the set consisting of
the empty word, and this defines a presentation for $\Gamma$.

	Let us call a word over $F$ \emph{trivial} if it corresponds
to the identity element of $\Gamma$.  Suppose that $w$ is a trivial
word, with
\begin{equation}
	w = \beta_1 \beta_2 \cdots \beta_n,
\end{equation}
where each $\beta_i$ is an element of $F$ or an inverse of an element
of $F$.  This leads to a path $\theta_0, \theta_1, \ldots, \theta_n$,
where $\theta_0$ is the identity element of $\Gamma$ and $\theta_j$ is
equal to $\beta_1 \beta_2 \cdots \beta_j$ when $j \ge 1$.  Because $w$
is a trivial word, $\theta_n$ is also the identity element in
$\Gamma$, which is to say that this path is a loop that begins and
ends at the identity element.

	Fix a finite set $R$ of relations, so that $F$ and $R$ give a
presentation for $\Gamma$.  Let $w$ be a trivial word over $F$ which
is also irreducible.  Define $A(w)$ to be the smallest nonnegative
integer $A$ for which there exist relations $r_1, r_2, \ldots, r_k$ in
$R$, integers $b_1, b_2, \ldots, b_k$, and words $u_1, u_2, \ldots,
u_k$ over $F$ such that the expression
\begin{equation}
  u_1 r_1^{b_1} u_1^{-1} u_2 r_2^{b_2} u_2^{-1} \cdots u_k r_k^{b_k} u_k^{-1}
\end{equation}
can be reduced to $w$ after cancellations as before,
\begin{equation}
	\sum_{j=1}^k L(u_j) \le A,
\end{equation}
and
\begin{equation}
	\sum_{j=1}^k |b_j| \, L(r_j)^2 \le A.
\end{equation}
Here if $z$ is a word over $F$ and $b$ is an integer, then $z^b$ is
defined in the obvious manner, by simply repeating $z$ $b$ times when
$b \ge 0$, or repeating $z^{-1}$ $-b$ times when $b < 0$.  A
representation of this type for $w$ necessarily exists, since $F$ and
$R$ give a presentation for $\Gamma$.

	The group $\Gamma$ is said to be \emph{hyperbolic} if there
is a nonnegative real number $C_0 \ge 0$ so that
\begin{equation}
	A(w) \le C_0 \, L(w)
\end{equation}
for all irreducible trivial words $w$.  The property of hyperbolicity
does not depend on the choice of finite presentation for $\Gamma$, and
in fact there are other definitions for which one only needs to assume
that $\Gamma$ is finitely generated, and the existence of a finite
presentation is then a consequence.  This characterization of
hyperbolicity is discussed in Section 2.3 of \cite{misha-hyp}.  Some
examples of hyperbolic groups are finitely-generated free groups and
the fundamental groups of compact connected Riemannian manifolds
without boundary and strictly negative curvature.  In particular,
this includes the fundamental group of a closed Riemann surface
with genus at least $2$.

	Let $M$ be a nonempty set.  A nonnegative real-valued function
$d(x,y)$ on the Cartesian product $M \times M$ is said to be a
\emph{quasimetric} if $d(x,y) = 0$ exactly when $x = y$, $d(x,y) =
d(y,x)$ for all $x, y \in M$, and
\begin{equation}
	d(x,z) \le C \Bigl(d(x,y) + d(y,z)\Bigr)
\end{equation}
for some positive real number $C$ and all $x, y, z \in M$.
If this last condition holds with $C = 1$, then $d(x,y)$ is
said to be a \emph{metric} on $M$.

	If $d(x,y)$ is a quasimetric on $M$ and $a$ is a positive real
number, then $d(x,y)^a$ is also a quasimetric on $M$.  If $d(x,y)$ is
a metric on $M$ and $a$ is a positive real number such that $a \le 1$,
then $d(x,y)^a$ is a metric on $M$ too.  These statements are not
difficult to verify.  There is a very nice result going in the other
direction, which states that if $d(x,y)$ is a quasimetric on $M$,
then there are positive real numbers $C'$, $\delta$ and a metric
$\rho(x,y)$ on $M$ such that
\begin{equation}
	C'^{-1} \, \rho(x,y)^\delta \le d(x,y) \le C' \, \rho(x,y)^\delta
\end{equation}
for all $x, y \in M$.  See \cite{MS1}.

	If $d(x,y)$ is a quasimetric on $M$, then one has many of the
same basic notions as for a metric, such as convergence of sequences,
open and closed sets, dense subsets, and so on.  For instance, it
makes sense to say that $M$ is separable with respect to a quasimetric
if it has a subset which is at most countable and also dense, and one
can define the topological dimension for $M$ as in \cite{HW}.  The
diameter of a subset can be defined in the usual manner using the
quasimetric, and this permits one to define the Hausdorff dimension of
a nonempty subset of $M$.  A famous result about metric spaces is that
the topological dimension is always less than or equal to the
Hausdorff dimension.  See Chapter VII of \cite{HW}.  This does not
work for quasimetrics in general, and it cannot possibly work.  For if
$(M, d(x,y))$ is a quasimetric space with Hausdorff dimension $s$ and
$a$ is a positive real number, then $(M, d(x,y)^a)$ has Hausdorff
dimension $s / a$, while the topological dimension of $(M, d(x,y)^a)$
is the same as that of $(M, d(x,y))$.

	Let $\Gamma$ be a finitely-presented group which is
hyperbolic.  Associated to $\Gamma$ is a space $\Sigma$ which is a
kind of ``space at infinity'' or ideal boundary of $\Gamma$,
consisting of equivalence classes of asymptotic directions in
$\Gamma$.  This space is a compact Hausdorff topological space of
finite dimension, as on p110-1 of \cite{misha-hyp}, and it contains a
copy of the Cantor set as soon as it has at least three elements.  If
$\Sigma$ has at most two elements, then $\Gamma$ is said to be
\emph{elementary}.  For a free group with at least two generators the
space at infinity is homeomorphic to a Cantor set, while ${\bf Z}$, a
free group with one generator, is elementary and has two points in the
space at infinity.  If $\Gamma$ is the fundamental group of a closed
Riemann surface of genus at least $2$, then $\Sigma$ is homeomorphic
to the unit circle in ${\bf R}^2$.  More generally, if $\Gamma$ is the
fundamental group of a compact $n$-dimensional Riemannian manifold
without boundary with strictly negative curvature, then $\Sigma$ is
homeomorphic to the unit sphere ${\bf S}^{n-1}$ in ${\bf R}^n$.

	Actually, the space at infinity is defined for any hyperbolic
metric space in \cite{misha-hyp}, and this can be specialized to a
hyperbolic group.  It is often preferable to work with metric spaces
which are ``geodesic'', in the sense that any pair of points can be
connected by a curve whose length is equal to the distance between the
two points.  It is often useful to think of a hyperbolic group as acting
on a geodesic hyperbolic metric space by isometries, and to use that
to study the space at infinity.

	It does not customarily seem to be said this way, but I think
it is fair to say that what are basically defined on the space at
infinity are quasimetrics, at least initially.  More precisely, it is
more like the logarithm of a quasimetric, or, in other words, there is
a one-parameter family of quasimetrics which are powers of each other.
A few years ago Gromov casually asked about approximating quasimetrics
by metric in the manner described before, and this is presumably the
reason.  In Section 7.2 of \cite{misha-hyp} one takes a different route,
in effect compactifying a geodesic hyperbolic metric space by looking
at modified measurements of lengths of curves which take densities
into account, densities that decay at infinity in a suitable manner.

	In nice situations, such as hyperbolic groups, and universal
coverings of compact Riemannian manifolds without boundary and
strictly negative curvature in particular, there are doubling
conditions on the space at infinity.  Compare with \cite{pierre1}.
There are also interesting measures around, as in \cite{coornaert}.

	A well-known result of Borel \cite{borel, raghunathan} says
that simply-connected symmetric spaces can be realized as the
universal covering of a compact manifold.  If the symmetric space is
of noncompact type and rank $1$, it has negative curvature, and thus
the fundamental group of the compact quotient, which is a uniform
lattice in the group of isometries of the symmetric space, is a
hyperbolic group.  If the symmetric space is a classical hyperbolic
space of dimension $n$, with constant negative curvature, then the
space at infinity can be identified with a Euclidean sphere of
dimension $n - 1$.  If the symmetric space is a complex hyperbolic
space of complex dimension $m$, then the space at infinity can be
identified topologically with a Euclidean sphere of real dimension $2m
- 1$, but the geometry corresponds to a sub-Riemannian space when $m
\ge 2$.  For other symmetric spaces of noncompact type and rank $1$,
one again obtains topological spheres of dimension $1$ less than the
real dimension of the symmetric space, with more complicated
sub-Riemannian structures.

\section{$p$-Adic numbers}
\label{p-adic numbers}
\setcounter{equation}{0}

	Let ${\bf Z}$ denote the integers, ${\bf Q}$ denote the
rational numbers, and let $|\cdot |$ denote the usual absolute value
function or modulus on the complex numbers ${\bf C}$.  On the rational
numbers there are other absolute value functions that one can
consider.  Namely, if $p$ is a prime number, define the \emph{$p$-adic
absolute value function} $|\cdot |_p$ on ${\bf Q}$ by $|x|_p = 0$ when
$x = 0$, $|x|_p = p^{-k}$ when $x = p^k m/n$, where $k$ is an integer
and $m$, $n$ are nonzero integers which are not divisible by $p$.  One
can check that
\begin{equation}
	|x y|_p = |x|_p \, |y|_p
\end{equation}
and
\begin{equation}
	|x + y|_p \le |x|_p + |y|_p
\end{equation}
for all $x, y \in {\bf Q}$, and in fact
\begin{equation}
	|x + y|_p \le \max(|x|_p, |y|_p)
\end{equation}
for all $x, y \in {\bf Q}$.

	Just as the usual absolute value function leads to the
distance function $|x - y|$, the $p$-adic absolute value function
leads to the $p$-adic distance function $|x - y|_p$ on ${\bf Q}$.
With respect to this distance function, the rationals are not complete
as a metric space, and one can complete the rationals to get a larger
space ${\bf Q}_p$.  This is analogous to obtaining the real numbers by
completing the rationals with respect to the standard absolute value
function.  By standard reasoning the arithmetic operations and
$p$-adic absolute value function extend from ${\bf Q}$ to ${\bf Q}_p$,
with much the same properties as before.  In this manner one gets the
field of $p$-adic numbers.  As a metric space, ${\bf Q}_p$ is complete
by construction, and one can also show that closed and bounded subsets
of ${\bf Q}_p$ are compact.  This is also similar to the real numbers.

	Note that the set ${\bf Z}$ of integers forms a bounded subset
of ${\bf Q}_p$, in contrast to being an unbounded subset of ${\bf R}$.
In fact, each integer has $p$-adic absolute value less than or equal
to $1$.  There are general results about absolute value functions on
fields to the effect that if the absolute values of integers are bounded,
then they are less than or equal to $1$, and the absolute value function
satisfies the ultrametric version of the triangle inequality.  See p28-9
of \cite{Gouvea}.  In this case the absolute value function is said to
be \emph{non-Archimedian}.  If the absolute values of integers are not
bounded, as in the case of the usual absolute value function, then the
absolute value function is said to be \emph{Archimedian}.

	A related point is that the set ${\bf Z}$ of integers is a
\emph{discrete} subset of the real numbers.  It has no limit points,
and in fact the distance between two distinct integers is always at
least $1$.  This is not the case in ${\bf Q}_p$, where ${\bf Z}$ is
bounded, and hence precompact.  Now consider ${\bf Z}[1/p]$, the set
of rational numbers of the form $p^k n$, where $k$ and $n$ are
integers.  As a subset of ${\bf R}$, this is unbounded, and it also
contains nontrivial sequences which converge to $0$.  Similarly, as a
subset of ${\bf Q}_p$, it is unbounded and contains nontrivial
sequences which converge to $0$.  As a subset of ${\bf Q}_l$ when $l
\ne p$, ${\bf Z}[1/p]$ is bounded and hence precompact again.  Using
the diagonal mapping $x \mapsto (x, x)$, one can view ${\bf Z}[1/p]$
as a subset of the Cartesian product ${\bf R} \times {\bf Q}_p$.  In
this product, ${\bf Z}[1/p]$ is discrete again.  Indeed, if $a = p^k
b$ is a nonzero element of ${\bf Z}[1/p]$, where $k$, $b$ are integers
and $b$ is not divisible by $p$, then either $|a|_p \ge 1$, or $|a|_p
\le 1$, in which case $k \ge 0$, and $|a| \ge 1$.

	Similarly, $SL_n({\bf Z})$, the group of $n \times n$
invertible matrices with entries in ${\bf Z}$ and determinant $1$, is
a discrete subgroup of $SL_n({\bf R})$, the analogously-defined group
of matrices with real entries.  One can define $SL_n({\bf Z}[1/p])$
and $SL_n({\bf Q}_p)$ in the same manner, and using the diagonal
embedding $x \mapsto (x, x)$ again, $SL_n({\bf Z}[1/p])$ becomes a
discrete subgroup of the Cartesian product $SL_n({\bf R}) \times
SL_n({\bf Q}_p)$.

	There are fancier versions of these things for making ${\bf
Q}$ discrete, using ``ad\`eles'', which involve $p$-adic numbers for all
primes $p$.  See \cite{Weil}.

	Now let us turn to some aspects of analysis.  With respect to
addition, ${\bf Q}_p$ is a locally compact abelian group, and thus has
a translation-invariant Haar measure, which is finite on compact sets,
strictly positive on nonempty open sets, and unique up to
multiplication by a positive real number.  As in \cite{Taibleson},
there is a rich Fourier analysis for real or complex-valued functions
on ${\bf Q}_p$, or ${\bf Q}_p^n$ when $n$ is a positive integer.

	Instead one can also be interested in ${\bf Q}_p$-valued
functions on ${\bf Q}_p$, or on a subset of ${\bf Q}_p$.  It is
especially interesting to consider functions defined by power series.
As is commonly mentioned, a basic difference between ${\bf Q}_p$ and
the real numbers is that an infinite series $\sum a_n$ converges if
and only if the sequence of terms $a_n$ tends to $0$ as $n$ tends to
infinity.  Indeed, the series converges if and only if the sequence of
partial sums forms a Cauchy sequence, and this implies that the terms
tend to $0$, just as in the case of real or complex numbers.  For
$p$-adic numbers, however, one can use the ultrametric version of the
triangle inequality to check that the partial sums form a Cauchy
sequence when the terms tend to $0$.  In particular, a power series
$\sum a_n \, x^n$ converges for some particular $x$ if and only if the
sequence of terms $a_n \, x_n$ tends to $0$, which is to say that
$|a_n|_p \, |x|_p^n$ tends to $0$ as a sequence of real numbers.

	Suppose that $\sum_{n = 0}^\infty a_n \, x^n$ is a power
series that converges for all $x$ in ${\bf Q}_p$, which is equivalent
to saying that $|a_n|_p r^n$ converges to $0$ as a sequence of real
numbers for all $r > 0$.  Thus we get a function $f(x)$ defined
on all of ${\bf Q}_p$, and we would like to make an analogy with
entire holomorphic functions of a single complex variable.  This
is somewhat like the situation of starting with a power series that
converges on all of ${\bf R}$, and deciding to interpret it as a 
function on the complex numbers instead.

	In fact, let us consider the simpler case of a power series
with only finitely many nonzero terms, which is to say a polynomial.
As in the case of complex numbers, it would be nice to be able to
factor polynomials.  The $p$-adic numbers ${\bf Q}_p$ are not
algebraically closed, and so in order to factor polynomials one 
can first pass to an algebraic closure.  It turns out that the $p$-adic
absolute value can be extended to the algebraic closure, while keeping
the basic properties of the absolute value.  See \cite{Cassels, Gouvea}.
The algebraic closure is not complete in the sense of metric spaces
with respect to the extended absolute value function, and one can
take a metric completion to get a larger field to which the absolute
values can be extended again.  A basic result is that this metric
completion is algebraically closed, so that one can stop here.
Let us write ${\bf C}_p$ for this new field, which is algebraically
closed and metrically complete.

	Once one goes to the algebraic closure, one can factor
polynomials.  On ${\bf C}_p$ one has this property and also one can
work with power series.  In particular, since the power series
$\sum_{n = 0}^\infty a_n \, x^n$ converges on all of ${\bf Q}_p$, it
also converges on all of ${\bf C}_p$, so that $f(x)$ can be extended
in a natural way to ${\bf C}_p$.  For that matter, one can start with
a power series that converges on all of ${\bf C}_p$, where the
coefficients are allowed to be in ${\bf C}_p$, and not just ${\bf Q}_p$.

	Under these conditions, the function $f$ can be written as
a product of an element of ${\bf C}_p$, factors which are equal to $x$,
and factors of the form $(1 - \lambda_j \, x)$, where the
$\lambda_j$'s are nonzero elements of ${\bf C}_p$.  In other words,
the factors of $x$ correspond to a zero of some order at the origin,
while the factors of the form $(1 - \lambda_j \, x)$ correspond
to zeros at the reciprocals of the $\lambda_j$'s.  If $f(x)$ is a
polynomial, then there are only finitely many factors, and this
statement is the same as saying that ${\bf C}_p$ is algebraically
closed.  In general, each zero of $f(x)$ is of finite order, and
there are only finitely many zeros within any ball of finite radius
in ${\bf C}_p$.  Thus the set of zeros is at most countable, and 
this condition permits one to show that the product of the factors
mentioned above converges when there are infinitely many factors.

	This representation theorem can be found on p113 of
\cite{Cassels} and on p209 of \cite{Gouvea}.  It is analogous to
classical results about entire holomorphic functions of a complex
variable, with some simplifications.  In the complex case, it is
necessary to make assumptions about the growth of an entire function
for many results, and the basic factors often need to be more
complicated in order to have convergence of the product.  See
\cite{Ahlfors, Veech} concerning entire holomorphic functions of
a complex variable.

\end{document}